\documentclass[a4paper,12pt]{article}
\usepackage{etex}
\usepackage[utf8x]{inputenc}
\usepackage{caption, subcaption}
\usepackage{arabtex}
\usepackage{hyperref}
\usepackage{relsize}
\usepackage{mathrsfs}
\usepackage[all]{xy}
\usepackage{amsmath,amssymb,amsfonts,amscd, graphicx, latexsym, verbatim, multirow, color}
\usepackage{epsfig}
\usepackage{geometry} 
\usepackage{adjustbox}

\newcommand{\F}{\mathcal F}

\newcommand{\Q}{\mathbb Q}
\newcommand{\Z}{\mathbb Z}

 \renewcommand{\R}{\mathbf R}

\newcommand{\psl}{\mathrm{PSL}_2(\Z)}
\newcommand{\pgl}{\mathrm{PGL}_2(\Z)}
\newcommand{\PGL}{\mathrm{PGL}}

\usepackage{color, ulem}

\newcommand\nt[1]{\textcolor{red}{{ #1}}}

\newcommand{\Jimm}{{\mathbf J}}

\newcommand{\sherh}[1]{\fboxsep=0pt\setlength{\fboxrule}{.5pt}
\begin{center}
   \fbox{\colorbox{white}{
         \begin{minipage}[t]{14.2cm}
            #1
         \end{minipage}
      }
   }
\end{center}}
\newcommand{\sherhh}[1]{\fboxsep=0pt\setlength{\fboxrule}{1pt}
\begin{center}
   \fbox{\colorbox{yellow}{
         \begin{minipage}[t]{13cm}
            #1
         \end{minipage}
      }
   }
\end{center}}

\newcommand{\sherhhh}[1]{\fboxsep=0pt\setlength{\fboxrule}{1pt}
\begin{center}
   \fbox{\colorbox{red}{
         \begin{minipage}[t]{13cm}
            #1
         \end{minipage}
      }
   }
\end{center}}

\renewcommand{\sherh}[1]{}\renewcommand{\sherhh}[1]{}\renewcommand{\sherhhh}[1]{}
\usepackage{forest}

\usepackage{tikz}
\usetikzlibrary{matrix,arrows}

\begin{document}

\title{A subtle symmetry of Lebesgue's measure}
\author{Muhammed Uluda\u{g}$^*$, Hakan Ayral\footnote{
{Department of Mathematics, Galatasaray University}
{\c{C}{\i}ra\u{g}an Cad. No. 36, 34357 Be\c{s}ikta\c{s}}
{\.{I}stanbul, Turkey}}}

\maketitle

\begin{abstract}
We represent the Lebesgue measure on the unit interval as a boundary measure of the Farey tree and show that this representation has a certain symmetry related to the tree automorphism induced by Dyer's outer automorphism of the group PGL(2,Z). Our approach gives rise to three new measures on the unit interval which are possibly of arithmetic significance.
\end{abstract}

\paragraph{\textsection1. Introduction.}\label{sec:introduction}
Lebesgue's measure on $\R$ enjoys two important invariance properties: $\lambda(x+A)=\lambda(A)$
and $\lambda(xA)=|x|\lambda(A)$ for $x\in \R$ and $A$ any measurable set. Our aim here is to represent Lebesgue's measure on the unit interval as a measure on the boundary of the Farey tree and show that it is symmetric under a natural operation which we shall introduce later in the paper.

Here, the Farey tree $\F$ is the first left-branch of Stern-Brocot's tree of rationals in $\Q\cap [0,1]$. If we specify left and right-turning probabilities at each vertex for a non-backtracking random walker starting at the top of  $\F$, we obtain a measure on the boundary $\partial\F$ inducing a natural Borel measure on the unit interval $[0,1]$. 
In this paper we consider the left/right turning probabilities inducing Lebesgue's measure on $[0,1]$. 

Automorphisms of $\F$ act by pre-composition of the left/right turning probabilities. 
The reflection of $\F$ around its symmetry axis is an automorphism $K$ acting by pre-composition, and the Lebesgue measure is invariant under it. 
On the other hand, every automorphism of $\F$ also act on the boundary $\partial\F$ and hence on $[0,1]$ via the continued fraction map. $K$ acts by sending the probability $r\in [0,1]$ to $1-r$, i.e. $K$ has a second action (by post-composition) on left/right turning probabilities. The Minkowski measure (see \textsection8) is the only measure left invariant by this latter action.

The subtle symmetry in question states that, in case of Lebesgue's measure, there is a second $\F$ -automorphism $\Jimm$ (jimm) acting both by pre- and post-composition and such that the two actions coincide. This automorphism is induced by Dyer's outer automorphism of the group $\PGL(2,\Z)$ (see \cite{dyer} and \cite{jimm}). The precise definition of this symmetry is given and proved in the final paragraph \textsection10, and a proper understanding of its statement should require all the preceding paragraphs. 

The paper is self-contained, except that the basic facts about the Stern-Brocot tree from \cite{orderings} and \cite{calk} are used freely in the text.

\paragraph{\textsection2. The Farey tree.} 
The {\it Farey tree}  $\F$ is the first left-branch of the the Stern-Brocot tree, whose vertices are labeled with the rational numbers in $(0,1)$.  It is generated from the rationals $0/1$ and $1/1$ by the well-known process of taking medians. 
If $ p /q < r/ s$ are consecutive fractions  in  $\F$, then $qr-ps=1$.

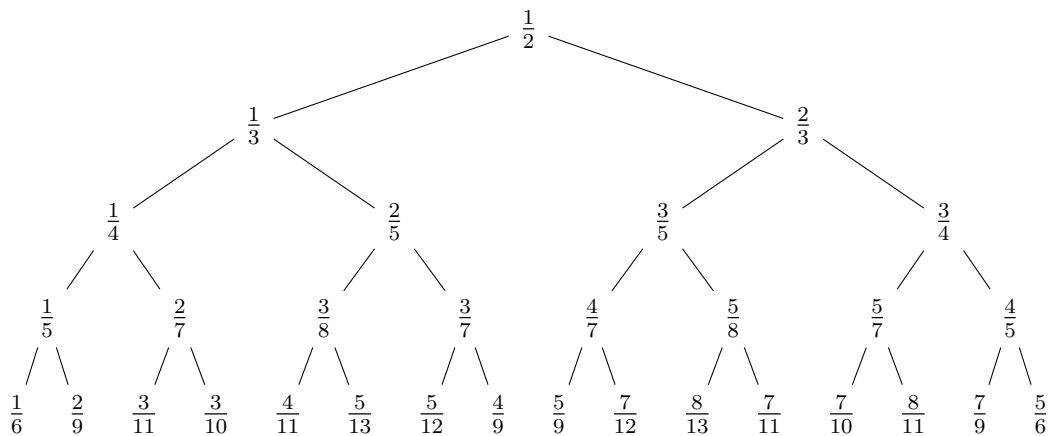
\begin{figure}[h]
\begin{forest}
  Stern Brocot/.style n args={5}{%
    content=$\frac{\number\numexpr#1+#3\relax}{\number\numexpr#2+#4\relax}$,
    if={#5>0}{
      append={[,Stern Brocot={#1}{#2}{#1+#3}{#2+#4}{#5-1}]},
      append={[,Stern Brocot={#1+#3}{#2+#4}{#3}{#4}{#5-1}]}
    }{}}
[,Stern Brocot={0}{1}{1}{1}{4}]
\end{forest}
\caption{\small The Farey tree}
\end{figure}
The vertex set of $\F$ is identified with  the set $\Q\cap (0,1)$. For every vertex $x\in \Q\cap (0,1)$, there is the corresponding {\it ocean} lying beneath that vertex in the tree, which we also label by $x$. The leftmost ocean is labeled $0/1$ and the rightmost one $1/1$. Every edge $I$ of the tree is an isthmus between the two seas $p/q$ and $r/s$, and we label this edge by the {\it Farey interval} $I=[p/q,r/s]\subset \R$ with $qr-ps=1$. Its length is 
$$
|I|= \frac p q - \frac r s =\frac{qr-ps}{qs}=\frac{1}{qs}.
$$ 
If (and only if) the path joining the vertex $u/v$ to the root passes through $I=[p/q,r/s]$, then $u/v\in I$. We may index the edges of $\F$ by the 
downward vertices of these edges, i.e. by the rationals in $\Q\cap (0,1)$. The index of $I=[p/q,r/s]$ is then the median $p/q \oplus r/s$.
For example, the index of $[0,1/(n-1)]$ is $1/n$.

\paragraph{\textsection3. A monoid structure on $\Q\cap (0,1)$.} 
Now we endow the set of vertices of $\F$ with a monoid structure. We identify each vertex of $\F$ (=element of $\Q\cap (0,1)$) by the path from the root to that vertex. To multiply vertices we concatenate the corresponding paths.
This is essentially the structure of the modular group $\mathrm{PSL}_2(\Z)$ transferred over to $\Q\cap (0,1)$.

To be more precise, consider the set 
$$
X:=\{(n_1, n_2, \dots n_k)\, |\, k, n_1, n_2,\dots, n_k=1,2,3,\dots\}
$$
For $x=(n_1, n_2, \dots n_k)\in X$ set $\|x\|:=n_1+n_2+ \dots +n_k-1$ and $\ell(x):=k$. 
The set $X$ is in bijection with $\Q\cap (0,1)$ via the continued fraction map
\begin{equation}
\theta: (n_1, n_2, \dots n_k)\in X \to [0,n_1, n_2, \dots n_k+1]=[0,n_1, n_2, \dots n_k,1]\in\Q\cap (0,1),
\end{equation}
where we denote as usual
$$
[0,n_1, n_2, \dots , n_k]=\cfrac{1}{n_1 + \cfrac{1}{\ddots \, + \cfrac{1}{n_{k}}}}.
$$
Every rational in $(0,1)$ is uniquely represented by a tuple in $X$. The two children of the node $[0,n_1,\dots, n_k+1]=[0,n_1,\dots, n_k,1]$ are $[0,n_1,\dots, n_k+2]$ and $[0,n_1,\dots, n_k,2]$.  In other words, the two children of the node $(n_1,\dots, n_k)$ are $(n_1,\dots, n_k+1)$ and $(n_1,\dots, n_k,1)$. 
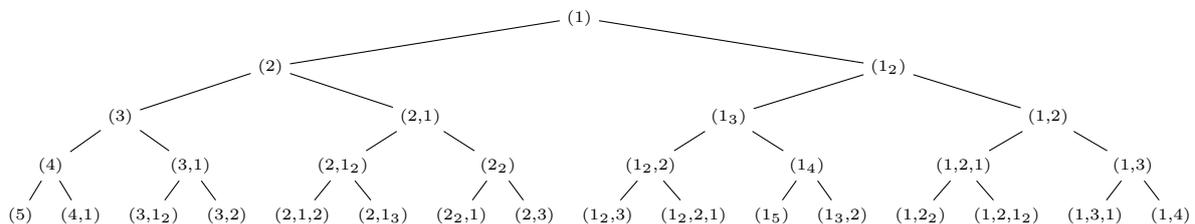
\begin{figure}[h]
\hspace{-.5cm}
{\tiny
\begin{forest}  T/.style n args={5}
[{(1)}
[{(2)}[{(3)}[{(4)}[{(5)}][{(4,1)}]][{(3,1)}[{(3,$1_2$)}][{(3,2)}]]]  [{(2,1)}[{(2,$1_2$)}[{(2,1,2)}][{(2,$1_3$)}]][{($2_2$)}[{($2_2$,1)}][{(2,3)}]]]] 
[{($1_2$)}[{($1_3$)}[{($1_2$,2)}[{($1_2$,3)}][{($1_2$,2,1)}]][{($1_4$)}[{($1_5$)}][{($1_3$,2)}]]]  [{(1,2)}[{(1,2,1)}[{(1,$2_2$)}][{(1,2,$1_2$)}]][{(1,3)}[{(1,3,1)}][{(1,4)}]]]] 
] ]\end{forest}}
\caption{\small The monoid $X$ ($1_k$ denotes the sequence ${1,1,\dots, 1}$ of length $k$). The map $\theta$ gives an isomorphism of the above tree with the Farey tree.}
\end{figure}
%
The concatenation operation on $X$ mentioned above is made precise as follows:
Let $x=(n_1, n_2, \dots n_k)$, $y=(m_1, m_2, \dots , m_l)$. 
We say that $x$ is a {\it right-child} if $k=\ell(x)$ is even and a {\it left-child}
otherwise.
\begin{equation}
x\star y:=\begin{cases}
(n_1, n_2, \dots n_{k-1}, n_k, m_1-1, m_2, m_3, \dots , m_l), &  x \mbox{ is a right-child} \\
(n_1, n_2, \dots n_{k-1}, n_k+m_1-1, m_2, \dots , m_l), & x\mbox{ is a left-child,} \\
\end{cases}
\end{equation}
where it is understood that $(\dots, m,n,k,\dots)=(\dots, m+k,\dots)$ if $n=0$.  
This is an associative operation and endows $X$ with the structure of a monoid.

\noindent
{\bf Examples.} One has
\begin{equation}
\underbrace{(1,1)\star(1,1)\star\dots\star (1,1)}_{n\mbox{ \scriptsize times }}=(1,n), \quad
\mbox{ and }\quad
\underbrace{(2)\star(2)\star\dots \star(2)}_{n\mbox{ \scriptsize times }}=(n+1).
\end{equation}
We transfer this operation  to an operation on $\Q\cap (0,1)$ via the bijection $\theta$. We denote this operation by $\star$ as well.

\noindent
{\bf Examples.} The two examples above becomes
\begin{equation}
\underbrace{{2 \over 3}\star {2 \over 3} \star\dots \star{2 \over 3}}_{n\mbox{ \scriptsize times }}={n+1\over n+2},
\quad\mbox{ and }\quad
\underbrace{{1 \over 3}\star {1 \over 3} \star\dots \star{1 \over 3}}_{n\mbox{ \scriptsize times }}={1\over n+2}
\end{equation}
The neutral element for this product is $e=(1)$ (which is a left-child). It corresponds to the top vertex $=1/2$ of the Farey tree. One has
$$
\forall x,y\in X, \quad \|x\star y\|=\|x\|+\|y\|.
$$
This monoid is freely generated by the elements
$L:=(2)$, $R:=(1,1)$, 
with 
$$
(n_1, n_2, \dots)=
L^{n_1-1}\star R^{n_2}\star \dots 
$$
An element $x\in X$ is a right (left) child if and only if its expansion in $\langle L, R\rangle$
ends with an $R$ ($L$); equivalently, $\ell(x)$ is even (odd).

We shall identify the set $\Q\cap (0,1)$, the vertex of set of $\F$, and the monoid $X$. The Farey intervals will  be indexed by the elements of $X$ via this identification and using the indexing given at the end of \textsection 2.
The Farey interval indexed by $(n_1,\dots, n_k)\in X$ will be denoted by $I(n_1,\dots, n_k)$.

\paragraph{\textsection4. The boundary of $\F$.}
There is the following passage from the boundary of $\F$ to the unit interval $[0,1]$. 
Recall that the boundary $\partial\F$ is the set of non-backtracking infinite paths (ends) based at the root vertex. 
Its natural topology is generated by the open sets $\mathcal O_I$, where  $\mathcal O_I$ is the set of ends through the edge $I$ of $\F$.
The space $\partial\F$ is a Cantor set. Furthermore, it has the ordering induced from the planar embedding of $\F$ and this ordering is compatible with its topology. An end can be encoded as an infinite sequence $\alpha:=(n_1, n_2, n_3, \dots )$ of positive integers, possibly terminating with $\infty$, and the map sending $\alpha$ to the continued fraction $[0,n_1,n_2, n_3, \dots ]$ is a continuous order-preserving surjection from $\partial\F$ onto the unit interval $[0,1]$. It is injective with the exception that the ends 
$[0,\dots, n_k+1, \infty]$ and $[0,\dots, n_k, 1,\infty]$ are sent to the same rational 
$[0,\dots, n_k+1]=[0,\dots, n_k, 1]$. The open set $\mathcal O_I$ is sent to the Farey interval 
$I=[p/q, r/s] \subset [0,1]$.

\paragraph{\textsection5. The automorphism $K$.} The monoid $(X, \star)$ has the automorphism $K$ exchanging its generators 
\begin{equation}
K(1,1):=(2), \quad K(2):=(1,1)
\end{equation}
We transfer $K$ to $\Q\cap (0,1)$ via $\theta$ and denote the resulting involution by $K$ again. 
If $x=\theta(n_1, n_2, \dots)=R^{n_1-1}\star L^{n_2} \dots,$ then we see that $K$ has a very nice form:
$$
K(n_1, n_2, \dots)=L^{n_1-1}\star R^{n_2}\star \dots = (1,n_1-1, n_2, \dots) \implies Kx=1-x.
$$
On $\F$, the involution $K$ is nothing but the reflection with respect to the symmetry axis.
As such, $K$ determines an automorphism of $\F$ as well. Moreover, it induces a homeomorphism of $\partial \F$ which induces the homeomorphism of the unit interval sending $x\in [0,1]$ to $1-x$.

\medskip
\noindent
{\bf Remark.} $K$ is related to the outer automorphism of the modular group $\psl$. It respects the topology of $\partial\F$ and reverses its ordering, inducing a homeomorphism of $[0,1]$.

\paragraph{\textsection6. The flip.}
The {\it flip} $\varphi$ is the involutive unary operation on $X$ defined as 
\begin{equation}
x=(n_1,n_2,\dots, n_k) \rightarrow \varphi(x)=(n_k,n_{k-1},\dots, n_1).
\end{equation}
Via the bijection $\theta$, we may transfer this involution to $\Q\cap (0,1)$ as
$$
\varphi([0,n_1,n_2,\dots, n_k])=[0,n_k-1, n_{k-1}, \dots n_2,n_1+1]
$$
where we assumed that $n_k>1$. 

\sherh{Expanding $x$ and $y$ in $\langle L, R \rangle$, it is readily seen that $\varphi(x\star y)=\varphi(y)\star \varphi(x)$. \nt{mistake}}

\medskip \noindent {\bf Examples.}
For $n>1$ one has $\varphi(1/n)=1/n$ and  $\varphi(n/(n+1))=\varphi([0,1,n])=[0,n-1,2]=2/(2n-1)$.

\medskip
\noindent
{\bf Remark.} The flip operation is related to the inversion in the modular group. It has no extension to $\partial\F$ nor to $[0,1]$ as it terribly violates the topology. In the appendix, we provide a Maple code to evaluate the flip directly on $\Q$.

\sherh{How does the flip etc behave with respect to the Farey sum operation?}

\sherh{Euler showed that the CFs $[ a_1, a_2, \dots a_n]$ and its reversal $[ a_n, a_{n-1}, \dots a_1]$ have the same numerators if the first and last terms in both CFs are non-zero. When we reverse the CF list of a fraction less than 1, the new fraction will have the same denominator if its CF also begins with a zero.}

%
%
\begin{figure}[h]
\centering
{\small
\begin{forest}  T/.style n args={5}
[$1 \over 2$
	[$1 \over 3$
		[$1 \over 4$
			[$1 \over 5$
				[$1 \over 6$]
				[$ 5\over 6 $]
			]
			[$ 4\over5 $
				[$ 5\over9 $]
				[$ 4\over9 $]
			]
		]
		[$ 3\over4 $
			[$ 4\over7 $
				[$ 4\over11 $]
				[$ 7\over11 $]
			]
			[$ 3\over7 $
				[$ 7\over10 $]
				[$ 3\over 10$]
			]
		]
	] 
	[$ 2\over 3$
		[$ 3\over 5$
			[$ 3\over8 $
				[$ 3\over 11$]
				[$ 8\over11 $]
			]
			[$ 5\over8 $
				[$ 8\over13 $]
				[$ 5\over13 $]
			]
		]
		[$2 \over5 $
			[$ 5\over7 $
				[$ 5\over12 $]
				[$ 7\over12 $]
			]
			[$ 2\over7 $
				[$ 7\over9 $]
				[$ 2\over 9$]
			]
		]
	] 
]\end{forest}
}
\caption{\small The flipped Farey tree. The two children of $p/q$ are $p/(p+q)$ and $q/(p+q)$. (this is a close relative of Calkin-Wilf's tree \cite{calk}).}
\end{figure}
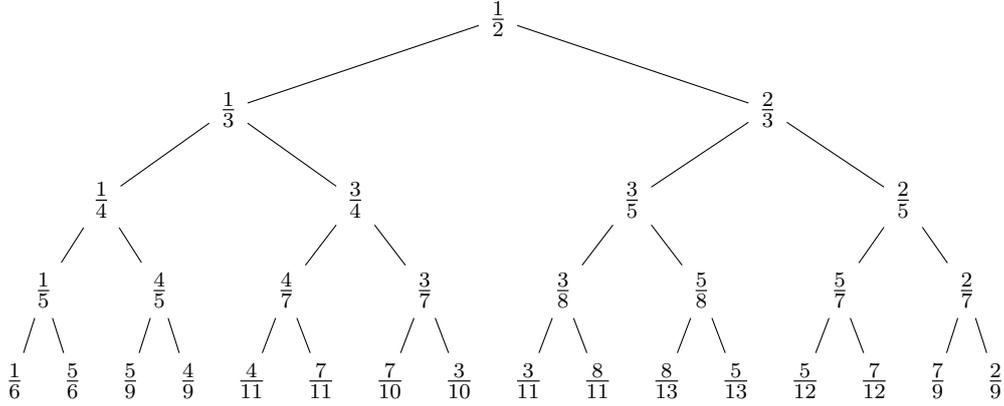
\paragraph{\textsection7. The involution Jimm.}
We define another unary operation on $X$ by
$$
x=(n_1,n_2,\dots, n_k) \rightarrow \Jimm(x)=(1_{n_1-1},2, 1_{n_2-2},2,1_{n_3-2},2,\dots, 1_{n_{k-1}-2},2, 1_{n_k-1}),
$$
where it is assumed that $k>1$ and the emerging $1_{-1}$'s are eliminated with the rule $[\dots m, 1_{-1},n,\dots]=[\dots m+n-1,\dots]$ and $1_0$ with the rule 
$[\dots m, 1_{0},n,\dots]=[\dots m,n,\dots]$, until all entries are $\geq 1$.  The latter rule includes the case $[1_0, m, \dots]=[m, \dots]$. These rules are applied once at a time. 
For $k=1$, set $\Jimm(x)=(1_{n_1})$.

\noindent
{\bf Examples.} We have\\
$\Jimm((1,1,1,1))=(1_{0},2, 1_{-1},2,1_{-1},2,1_0)=(3,1_{-1},2,1_0)=(4,1_0)=(4)$, \\
$\Jimm((2,2,2,2))=(1_{1},2, 1_{0},2,1_{0},2,1_1)=(1,2,2,1_{0},2,1)=(1,2,2,2,1)$,\\
$\Jimm((1,2,2,2,1))=(1_{0},2, 1_{0},2,1_{0},2,1_0,2,1_0)=(2,2,2,2)$.\\

Loosely speaking, $\Jimm$ sends the zig-zag segments in the path of $x$ to straight segments and vice versa.
Observe that $\|\Jimm(x)\|=\|x\|$ for any $x\in X$, i.e. $\Jimm$ preserves the depth on $\F$, but not always the length $\ell(x)$;  at a given depth it sends shorter tuples to longer tuples. 

Every automorphism of the Farey tree determines a unique bijection of $\Q\cap (0,1)$, i.e. there is a map
$Aut(\F)\to Sym(\Q\cap (0,1))$.
Conversely, a permutation of the set $\Q\cap (0,1)$ determines an automorphism of $\F$ if it preserves siblings.
As it obviously preserves siblings, $\Jimm$ defines an automorphism of $\F$.

Via $\theta$, we may transfer $\Jimm$ to $\Q\cap (0,1)$ as
$$
\Jimm([0,n_1,n_2,\dots, n_k])=[0,1_{n_1-1},2, 1_{n_2-2},2,1_{n_3-2},2,\dots, 1_{n_{k-1}-2},2, 1_{n_k-1}],
$$
where it is assumed that $n_k>1$. 

\noindent
{\bf Examples.} We have
$\Jimm((1))=(1) \iff \Jimm(1/2)=1/2$, 
$\Jimm((1,1))=(2)\iff \Jimm(1/3)=2/3$ and
$\Jimm((2))=(1,1) \iff \Jimm(2/3)=1/3.$

\bigskip
Observe  that, for $x=[0,n_1,n_2,\dots, n_k]$, one has 
\begin{eqnarray*}
\frac{x}{1+x}=[0,n_1+1,n_2,\dots, n_k]
\implies 
\Jimm\left(\frac{x}{1+x}\right)=[0,1_{n_1},2, 1_{n_2-2},\dots,2, 1_{n_k-1}]\\
=[0,1,1_{n_1-1},2, 1_{n_2-2},\dots,2, 1_{n_k-1}]=
\frac{1}{1+{\Jimm(x)}}
\end{eqnarray*}
Hence, we have the functional equations
\begin{eqnarray}\label{fe}
\Jimm\left(\frac{1}{1+x}\right)=\frac{\Jimm(x)}{1+\Jimm(x)}, \quad \Jimm\left(\frac{x}{1+x}\right)=\frac{1}{1+\Jimm(x)}.
\end{eqnarray}

\medskip
The following simple fact is a direct consequence of definitions:

\medskip
\noindent
{\bf Observation.} The involutions $\varphi$ and $K$  commutes with  $\Jimm$.

In the appendix, we provide a Maple code to evaluate $\Jimm$ directly on $\Q$.

\medskip
\noindent
{\bf Remark 1.} $\Jimm$ extends in a natural manner
to $[1,\infty)\cap \Q$ via $\Jimm(1/x):=1/\Jimm(x)$ and to $(-\infty, 0]\cap \Q$ via 
$\Jimm(-x):=-1/\Jimm(x)$. Applying $\Jimm$ to the Calkin-Wilf sequence gives another ``{twisted}" enumeration of $\Q^+$:
$$
\frac{1}{1}, \frac{1}{2}, \frac{2}{1}, \frac{2}{3}, \frac{3}{1}, \frac{1}{3}, \frac{3}{2}, \frac{3}{5}, \frac{5}{2}, \frac{1}{4}, \frac{4}{3}, 
\frac{3}{4}, \frac{4}{1}, \frac{2}{5}, \frac{5}{3}, \frac{5}{8}, \frac{8}{3}, \frac{2}{7}, \frac{7}{5}, \frac{4}{5}, \frac{5}{1}, \frac{3}{7}, 
\frac{7}{4}, \frac{4}{7}, \frac{7}{3}, \frac{1}{5}, \frac{5}{4}, \frac{5}{7}, \frac{7}{2}, \frac{3}{8},\dots
$$
\sherh{
$$
\frac{1}{1}, \frac{1}{2}, \frac{2}{1}, \frac{2}{3}, \frac{3}{1}, \frac{1}{3}, \frac{3}{2}, \frac{3}{5}, \frac{5}{2}, \frac{1}{4}, \frac{4}{3}, \frac{3}{4}, \frac{4}, \frac{2}{5}, \frac{5}{3}, \frac{5}{8}, \frac{8}{3}, \frac{2}{7}, \frac{7}{5}, \frac{4}{5}, \frac{5}, \frac{3}{7}, \frac{7}{4}, \frac{4}{7}, \frac{7}{3}, \frac{1}{5}, \frac{5}{4}, \frac{5}{7}, \frac{7}{2}, \frac{3}{8}, \frac{8}{5}, \frac{8}{13}, \frac{13}{5}, \frac{3}{11}, \frac{11}{8}, \frac{7}{9}, \frac{9}{2}, \frac{5}{12}, \frac{12}{7}, \frac{5}{9}, \frac{9}{4}, \frac{1}{6}, \frac{6}{5}, \frac{7}{10}, \frac{10}{3}, \frac{4}{11}, \frac{11}{7}, \frac{7}{11}, \frac{11}{4}, \frac{3}{10}, \frac{10}{7}, \frac{5}{6}, \frac{6}{1}, \frac{4}{9}, \frac{9}{5}, \frac{7}{12}, \frac{12}{5}, \frac{2}{9}, \frac{9}{7}, \frac{8}{11}, \frac{11}{3}, \frac{5}{13}, \frac{13}{8}, \frac{13}{21}, \frac{21}{8}, \frac{5}{18}, \frac{18}{13}, \frac{11}{14}, \frac{14}{3}, \frac{8}{19}, \frac{19}{11}, \frac{9}{16}, \frac{16}{7}, \frac{2}{11}, \frac{11}{9}, \frac{12}{17}, \frac{17}{5}, \frac{7}{19}, \frac{19}{12}, \frac{9}{14}, \frac{14}{5}, \frac{4}{13}, \frac{13}{9}, \frac{6}{7}, \frac{7}{1}, \frac{5}{11}, \frac{11}{6}, \frac{10}{17}, \frac{17}{7}, \frac{3}{13}, \frac{13}{10}, \frac{11}{15}, \frac{15}{4}, \frac{7}{18}, \frac{18}{11}, \frac{11}{18}, \frac{18}{7}, \frac{4}{15}, \frac{15}{11}, \frac{10}{13}
$$}

\medskip
\noindent
{\bf Remark 2.} In fact, $\Jimm$ involutive and it is related to the outer automorphism of $\pgl$. It respects the topology of $\partial\F$ though not its ordering. For every $x\in \R\setminus \Q$, the limit $\lim_{q\to x}\Jimm(q)$ exists and the limiting involution 
$\Jimm_\R$ is continuous at irrationals, jumps at rationals and is a.e. differentiable with a derivative vanishing a.e.. Equations(\ref{fe}) hold for $\Jimm_\R$ as well. Furthermore, $\Jimm_\R$ preserves the quadratic irrationals set-wise, commuting with the quadratic conjugation. It also respects the $\psl$-action on $\R$. For proofs of these facts and more details about $\Jimm$ and $\Jimm_\R$, see \cite{jimm} (beware the use of a different notation for $\Jimm$ in that paper). 

 %
 %
%
\begin{figure}[h]
\centering
{\small
\begin{forest}  T/.style n args={5}
 [$ 1 \over 2$ [$ 2 \over 3$ [$ 3 \over 5$ [$ 5 \over 8$ [$ 8 \over 13 $]  [$ 7 \over 11 $]  ]  [$ 4 \over 7$ [$ 5 \over 9 $]  [$ 7 \over 12 $]  ]  ]  [$ 3 \over 4$ [$ 4 \over 5$ [$ 7 \over 9 $]  [$ 5 \over 6 $]  ]  [$ 5 \over 7$ [$ 7 \over 10 $]  [$ 8 \over 11 $]  ]  ]  ]  [$ 1 \over 3$ [$ 1 \over 4$ [$ 2 \over 7$ [$ 3 \over 11 $]  [$ 3 \over 10 $]  ]  [$ 1 \over 5$ [$ 1 \over 6 $]  [$ 2 \over 9 $]  ]  ]  [$ 2 \over 5$ [$ 3 \over 7$ [$ 5 \over 12 $]  [$ 4 \over 9 $]  ]  [$ 3 \over 8$ [$ 4 \over 11 $]  [$ 5 \over 13 $]  ]  ]  ]  ] 
\end{forest}
}
\caption{\small The Jimm-transform of the Farey tree.} 
\end{figure}
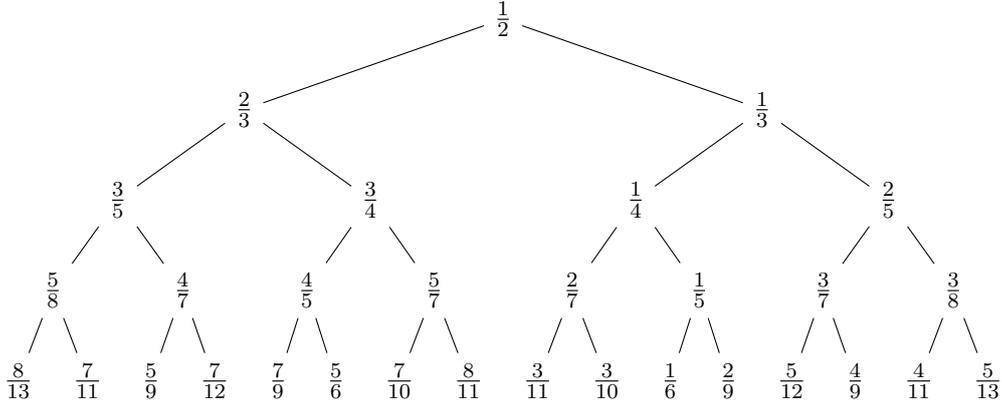

\sherh{
Conjugation provides an alternative enumeration of $\Q\cap(0,1)$.}

\paragraph{\textsection8. Measures on $\partial \F$.}
Imagine a random walker $W$ on $\F$ starting at the root and advancing towards the absolute without 
backtracking.
For every vertex $p$, we are given the probability $\pi(p)$ of arriving to that vertex from its parent. 
Hence we have a function $\pi: \Q\cap (0,1)\to [0,1]$ with the property: if $p$ and $q$ are siblings then $\pi(p)+\pi(q)=1$. 
For the root vertex, set $\pi(1):=1$. A function $\pi$ with these properties will be called a transition function. 

Let $\mathcal P$ be the set of all transition functions $\pi: \Q\cap (0,1)\to [0,1]$ on the set of vertices of $\F$.
Then $Aut(\F)$ acts on $\mathcal P$ by pre-composition, i.e. $\alpha\in Aut(\F)$ sending $\pi$ to
$\pi\circ \alpha$. On the other hand, if $\alpha$ commutes with $K\in Aut(\F)$ (there are many such), then $\alpha \circ \pi$ 
is also in $\mathcal P$ (provided $\alpha$ is well-defined on the image set of $\pi$).

The probability that the walker ends up in the 
Farey interval $I=[p/q,r/s]$ is the product of probabilities of choices he makes to go from the root to the vertex $p/q \oplus r/s$.
Since the boundary topology is generated by the Farey intervals, this puts a probability measure on the Borel algebra of $\partial\F$ and since the map $\partial\F\to [0,1]$ is measurable with respect to Borel algebras, we obtain a probability measure $\mu_\pi$ on $[0,1]$. Here, we assumed that $\mu$ has no point masses as it allows us to be sloppy about the endpoints of the Farey intervals. This is not an essential assumption, however.

Let ${\mathbf F}_\pi$ be the cumulative distribution function (c.d.f.) of $\mu_\pi$.
Since $\partial\F$ is an ordered space, inducing on $[0,1]$ its linear ordering, we will take the liberty to consider 
${\mathbf F}_\pi$ simultaneously as a c.d.f. on $\partial\F$ and on $[0,1]$.

Consider the ``{flipped Farey map}"
\begin{eqnarray}
 T_{\varphi F}: (n_1, n_2, \dots, n_{k-1}, n_k) \in X \to
 (n_1, n_2, \dots, n_{k-1}, n_k-1) \in X,
\end{eqnarray}
 where a zero at the end of a tuple is ignored. Then $T_{\varphi F}x$ is the parent of $x$ in $\F$.  
The naming of $ T_{\varphi F}$ is due to the fact that $ T_{F}:=\varphi  T_{\varphi F} \varphi$  is the usual {\it Farey map}
\begin{eqnarray}\label{farey1}
T_{F}: (n_1, n_2, \dots, n_{k-1}, n_k) \in X \to
 (n_1-1, n_2, \dots, n_{k-1}, n_k) \in X,
\end{eqnarray}
where a zero at the beginning of a tuple is ignored. 
Indeed, if we express $T_{F}$ in terms of continued fractions, we get the usual Farey map (\cite{orderings}):
\begin{eqnarray}\label{farey}
 T_{F}(x):=
 \begin{cases}
x/(1-x), & x<1/2\\
(1-x)/x, &x\geq 1/2.
\end{cases}
\end{eqnarray}
Now consider the interval $I(n_1, \dots, n_k)$. Then 
$$
\mu_\pi(I(n_1, \dots, n_k))=\prod_{i=0}^{d-1}\pi(T_{\varphi F}^{i}(n_1, n_2, \dots, n_k)),
 $$
 where $d=\|n_1, \dots, n_k\|+1=n_1+\dots+n_k$. In a more compact form, we may write
 \begin{eqnarray}
x\in X\implies \mu_\pi(x)=\prod_{i=0}^{\|x\|}\pi(T_{\varphi F}^{i}(x)).
\end{eqnarray}
 Suppose now that $x\in (0,1)$, $x=[0,n_1,n_2,\dots]$. Then for the c.d.f. of $\mu_\pi$ one has
\begin{eqnarray}
\mathbf F_\pi(x)=\mu_\pi([0,x])=\sum_{k=1}^\infty (-1)^k\mu_\pi\{W\in I(n_1, n_2, \dots, n_k)\}\\
=\sum_{k=1}^\infty (-1)^{1+k}\prod_{i=0}^{d-1}\pi(T_{\varphi F}^{i}(n_1, n_2, \dots, n_k)).
\end{eqnarray}

Note that every $\alpha\in Aut(\F)$ commutes with $T_{\varphi F}$ (but not always with $T_{F}$). This implies, for the c.d.f. of the measure defined by $\pi \circ \alpha$, 
$$
\mathbf F_{\pi \circ \alpha}(x)=\mathbf F_\pi(\alpha x).
$$

\noindent
{\bf Example: Minkowski's measure and Denjoy's measures.}
When $\pi\equiv 1/2$ (except the root vertex), then $\mu_\pi$ is ``{Minkowski's measure}" as $\mathbf{F}_\pi$ is nothing but the Minkowski question mark function \cite{alkauskas}, \cite{linas}. This measure is invariant under the full $Aut(\F)$-action and is the only measure with this property. One has 
 $$
\mu_\pi\{W\in I(n_1, \dots, n_k)\}=2^{1-n_1-n_2\dots-n_k}.
 $$
Its  c.d.f.  is
\begin{eqnarray}
\mathbf F_\pi(x)=\sum_{k=1}^\infty (-1)^{1+k}2^{1-n_1-n_2\dots-n_k}.
\end{eqnarray}
When $\pi$ assumes a constant value $a$ on left-children
(and thus the constant value $1-a$ on right-children), then the resulting measure generalising Minkowski's is called {\it Denjoy's measure}.
\paragraph{\textsection9. Lebesgue's measure.}
Denote this measure by $\lambda$ and its  $\pi$-function by $\pi_\lambda$.  
Let $x=(n_1,\dots, n_k)\in X$ with  $\theta(x)=[0,n_1,\dots,n_{k-1}, n_k+1]$ and let $I(x)$ be the corresponding Farey interval. 
Then 
$\mu_\pi(I(x))$ equals the length of $I(x)$, which is
$$
\lambda( I(x))=\prod_{i=0}^{\|x\|}\pi_\lambda(T_{\varphi F}^{i}(x)).
 $$
 Denote $[0,n_1,\dots, n_k]=:p_k/q_k$. Then $p_k$ and $q_k$ satisfy the recursions 
 $$
 p_k=n_kp_{k-1}+p_{k-2},\quad q_k=n_kq_{k-1}+q_{k-2}.
 $$
The endpoints of $I(x)$ are $[0,n_1,\dots, n_{k-1}]$ and $[0,n_1,\dots, n_{k-1}, n_{k}]$, so
\begin{eqnarray}\label{eqA}
\lambda(I(x))=
\left| \frac{n_kp_{k-1}+p_{k-2}}{n_kq_{k-1}+q_{k-2}}-\frac{p_{k-1}}{q_{k-1}}\right|=
\frac{1}{q_{k-1}q_{k}}=\\
\frac{1}{\langle n_1,  \dots, n_{k-1}\rangle \langle n_1,  \dots, n_{k-1}, n_{k}\rangle }.
\end{eqnarray}
Now let $y$ be the parent of $x$. 
Then $y=T_{\varphi F} x=(n_1, \dots, n_k-1)$ (zeros at the end are ignored) and  $\theta(y)=[0,n_1,\dots,n_{k-1}, n_k]$. Since the end points of $I(y)$ are $[0,n_1,\dots, n_{k}-1]$ and $[0,n_1,\dots, n_{k-1}]$,
\begin{equation}\label{eqB}
\lambda(I(y))=\left| \frac{(n_k-1)p_{k-1}+p_{k-2}}{(n_k-1)q_{k-1}+q_{k-2}}-\frac{p_{k-1}}{q_{k-1}}\right|
=\frac{1}{\langle n_1,  \dots, n_{k-1}\rangle \langle n_1, \dots, n_{k-1}, n_k-1 \rangle }.
\end{equation}
 From (\ref{eqA}) and (\ref{eqB}) we obtain
 $$
 \pi_\lambda(x)=
 \frac{\lambda(I(x))}{\lambda(I(y))}=
 \frac{\langle n_1,  \dots,n_{k-1},n_k-1\rangle }{\langle n_1, \dots,  n_{k}\rangle}=
 \frac{\langle n_1, \dots, n_k\rangle-\langle n_1, \dots, n_{k-1} \rangle}{\langle n_1, \dots, n_k \rangle}
 $$
 $$
= 1-
 \frac{\langle n_1, \dots, n_{k-1} \rangle}{\langle n_1, \dots, n_{k-1},n_k\rangle}=1-[0,n_{k}, n_{k-1},\dots, n_1].
 $$
 Hence, with our usual convention $\dots, m,0,k, {\dots} =\dots m+k\dots$ we may write
\begin{equation}\label{pilambda}
  \pi_\lambda(x)=[0,1,n_{k}-1, n_{k-1},\dots, n_1]
\end{equation}
We may express this in terms of $r\in \Q\cap (0,1)$ as: 
\begin{equation}\label{pilambda2}
\pi_\lambda[0,n_1,\dots, n_k]=1-[0, n_k-1,n_{k-1},\dots,  n_1]
\iff
\pi_\lambda(r)=1-\varphi(T_{F}r)=K\varphi T_{F}(r)
\end{equation}
Alternatively, we may express this as a map $X\to X$ as 
$$
{\pi_\lambda(n_1, \dots, n_k)=(1,n_k-1, n_{k-1}, \dots,  n_1-1),}
$$
where the zeros at the end are ignored.
In this description, the two-to-one nature of $\pi_\lambda$ becomes evident:
$$
\pi_\lambda(n_1, \dots, n_k)=\pi_\lambda(1,n_1-1, \dots, n_k) \iff \pi_\lambda(x)=\pi_\lambda(Kx)
$$
This is simply the $K$-invariance of Lebesgue's measure.
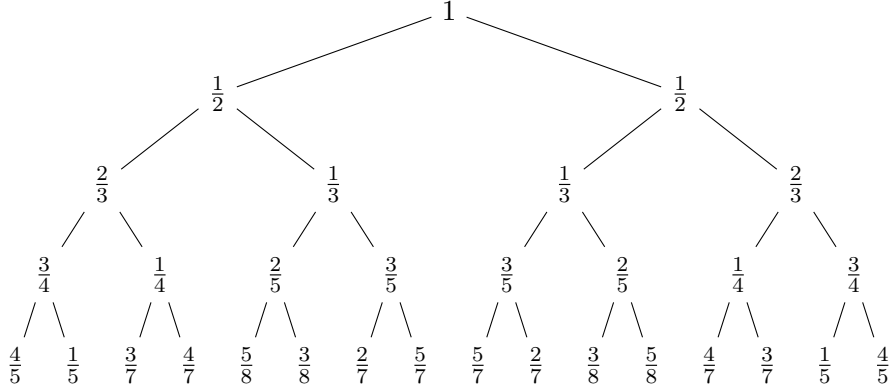
\begin{figure}[h]
\centering
{\small
\begin{forest}  T/.style n args={5}
[$1$
	[$1 \over 2$
		[$2 \over 3$
			[$3 \over 4$
				[$4 \over 5$
				]
				[$ 1\over5 $
				]
			]
			[$ 1\over4 $
				[$ 3\over7 $
				]
				[$ 4\over7 $
				]
			]
		] 
		[$ 1\over 3$
			[$ 2\over 5$
				[$ 5\over8 $
				]
				[$ 3\over8 $
				]
			]
			[$3 \over5 $
				[$ 2\over7 $
				]
				[$ 5\over7 $
				]
			]
		] 
	]
	[$1 \over 2$
		[$1 \over 3$
			[$3 \over 5$
				[$5 \over 7$
				]
				[$ 2\over7 $
				]
			]
			[$ 2\over5 $
				[$ 3\over8 $
				]
				[$ 5\over8 $
				]
			]
		] 
		[$ 2\over 3$
			[$ 1\over 4$
				[$ 4\over7 $
				]
				[$ 3\over7 $
				]
			]
			[$3 \over4 $
				[$ 1\over5 $
				]
				[$ 4\over5 $
				]
			]
		] 
	]
]
\end{forest}}
\caption{\small The ``{Lebesgue tree}" $\mathcal L$. This is obtained by replacing each node of the Farey tree, by the probability of arriving to that node from its parent. The two nephews (offsprings of a sibling) of $p/q$ are $q/(p+q)$ and $q/(p+q)$ in 
$\mathcal L$. Multiplying the numbers in a lineage, yields always a reciprocal integer.}
\end{figure}
%
\paragraph{\textsection10. The symmetry.}
Here is the symmetry of Lebesgue's measure promised in the title, enshrouded deeply in the Farey tree:
$$
\pi_\lambda\Jimm(x)=\Jimm\pi_\lambda(x)
$$
This is because, as shown in (\ref{pilambda2}), $\pi_\lambda=K\varphi T_{F}$ and $\Jimm$ commutes both with  $K$, $\varphi$ and $T_{F}$. The commutativity with $K$ and $\varphi$ was observed in \textsection7.
It remains to see that $\Jimm$ commutes with the Farey map $T_{F}$, which is easily observed from the description of $T_{F}$ in (\ref{farey1}).

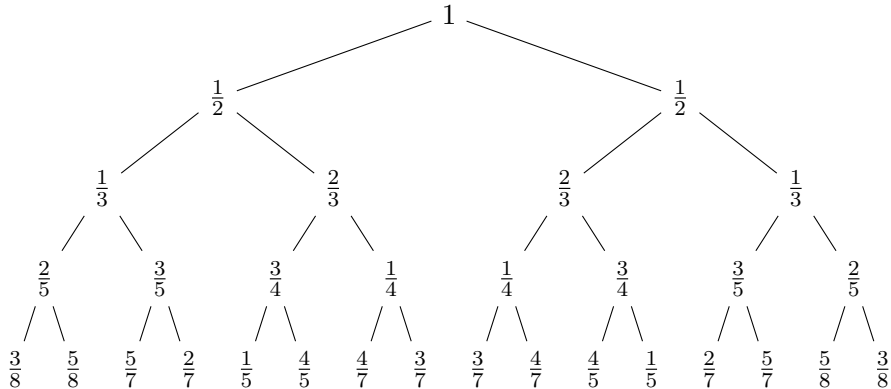
\begin{figure}[h]
\centering
{\small
\begin{forest}  T/.style n args={5}
[$1$
	[$1 \over 2$
		[$1 \over 3$
			[$2 \over 5$
				[$3 \over 8$			]
				[$5 \over 8$
				]
			]
			[$ 3\over5 $
				[$ 5\over7 $
				]
				[$ 2\over7 $				]
			]
		] 
		[$ 2\over 3$
			[$ 3\over 4$
				[$ 1\over5 $
				]
				[$ 4\over5 $
				]
			]
			[$1 \over4 $
				[$ 4\over7 $
				]
				[$ 3\over7 $
				]
			]
		] 
	]
	[$1 \over 2$
		[$2 \over 3$
			[$1 \over 4$
				[$3 \over 7$			]
				[$4 \over 7$
				]
			]
			[$ 3\over4 $
				[$ 4\over5 $
				]
				[$ 1\over5 $				]
			]
		] 
		[$ 1\over 3$
			[$ 3\over 5$
				[$ 2\over7 $
				]
				[$ 5\over7 $
				]
			]
			[$2 \over5 $
				[$ 5\over8 $
				]
				[$ 3\over8 $
				]
			]
		] 
	]
]
\end{forest}}
\caption{\small The jimm-transform of Lebesgue's tree.  The two nieces (offsprings of a sibling) of $p/q$ are $q/(p+q)$ and $q/(p+q)$. Multiplying the numbers in a lineage, yields always a reciprocal integer.}
\end{figure}
%

\paragraph{\textsection11. Conclusion.}
We are currently trying to understand how this symmetry manifests itself on the superficial level, i.e. on the arithmetic-related analysis and dynamics of the unit interval.  There are many questions pertaining to the measures induced by the transition functions
 $\pi(x):=K\pi_\lambda(x)$,  $\pi(x):=\Jimm\pi_\lambda(x)=\pi_\lambda\Jimm(x)$
and $\pi(x):=K\Jimm\pi_\lambda(x)=\Jimm K\pi_\lambda(x)$. These are, in a sense, basic deformations of Lebesgue's measure.
Their c.d.f.'s are depicted below (Figs 7-9); see our forthcoming paper \cite{deformations} for more details.

\begin{figure}[h]
    \centering
	\includegraphics[width=\textwidth]{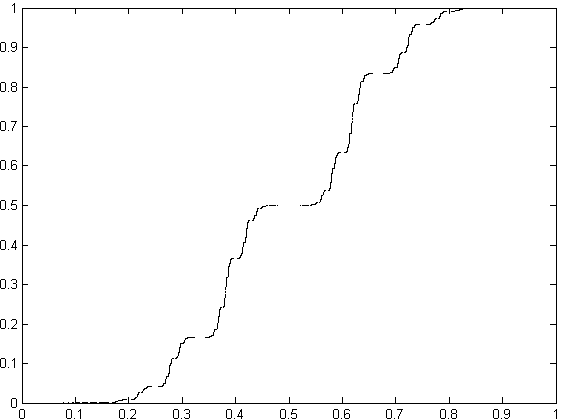}
	\caption{ c.d.f. of $K\pi_\lambda(x)$}
    \end{figure}

\begin{figure}[h]
    \centering
	\includegraphics[width=\textwidth]{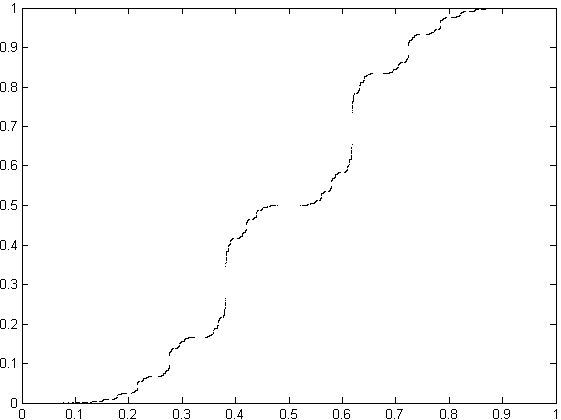}
	\caption{ c.d.f. of $\Jimm\pi_\lambda(x)$}
    \end{figure}

\begin{figure}[h]
    \centering
	\includegraphics[width=\textwidth]{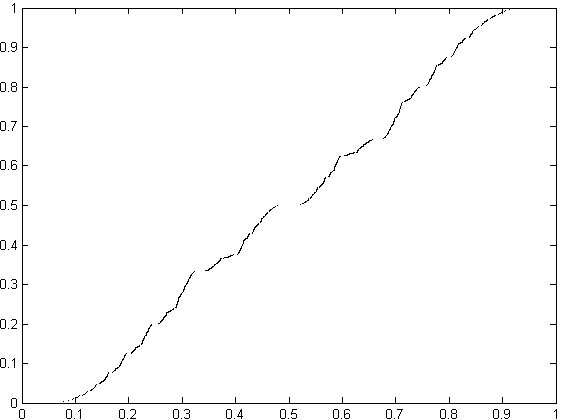}
	\caption{ c.d.f. of $K\Jimm \pi_\lambda(x)$}
    \end{figure}

\bigskip\noindent
{\bf Acknowledgements.}
This research is funded by a Galatasaray University research grant and  T\"UB\.ITAK grant 115F412. 
The second named author was also funded by the T\"UB\.ITAK grant 113R017.
We are grateful to an anonymous referee for shortening the proof of the expression (\ref{pilambda}).


\paragraph{Appendix: Maple code to evaluate Jimm and the Flip}
\begin{verbatim}
>with(numtheory)
>jimm := proc (q) local M, T, U, i, x; 
	T := matrix([[1, 1], [1, 0]]); 
	U := matrix([[0, 1], [1, 0]]); 
	M := matrix([[1, 0], [0, 1]]); 
	x := cfrac(q, quotients); 
	if x[1] = 0 then for i from 2 to nops(x) do 
		M := evalm(`&*`(`&*`(M, T^x[i]), U)) end do; 
	return M[2, 2]/M[1, 2] else for i to nops(x) do 
		M := evalm(`&*`(`&*`(M, T^x[i]), U)) end do; 
	return M[1, 2]/M[2, 2] end if 
end proc;

>flip := proc (q) local x, y, n; 
y := [0]; 
x := cfrac(q, quotients); 
n := nops(x); 
y := [op(y), x[n]-1]; 
for i to n-3 do y := [op(y), x[n-i]] end do;
 y := [op(y), x[2]+1]; 
 return cfrac(y) 
 end proc;
\end{verbatim}


\begin{thebibliography}{99}

\bibitem{alkauskas}
G. Alkauskas, 
{\it The Minkowski question mark function: explicit series for the dyadic period function and moments }
Math. Comp. 79 (269) (2010), 383-418; Addenda and corrigenda, Math. Comp. 80 (276) (2011), p.2445-2454.

\bibitem{Bruce} 
Reznick, Bruce. {\it  Regularity properties of the Stern enumeration of the rationals.} 
Journal of integer sequences 11.2 (2008): 3.
 

\bibitem{calk}
B. Bates, M. Bunder and K. Tognetti,
{\it Linking the Calkin--Wilf and Stern--Brocot trees}
European Journal of Combinatorics, 31 no.7 (2010) p.1637--1661.

\bibitem{orderings} 
C. Bonanno and S. Isola,
{\it Orderings of the rationals and dynamical systems},
arXiv preprint arXiv:0805.2178 (2008).

\bibitem{dyer} 
J. L. Dyer,
{\it Automorphism sequences of integer unimodular groups},
Illinois Journal of Mathematics, vol.22, no.1 (1978) p.1--30.


\bibitem{jimm} 
A.M. Uludağ and H.Ayral,
{\it Jimm, a Fundamental Involution}
arXiv preprint arXiv:1501.03787 (2015).
  
 \bibitem{deformations} 
A.M. Uludağ and H.Ayral,
{\it Some deformations of Lebesgue’s measure on the boundary of the Farey tree,}
to appear.
  
\bibitem{linas} 
L. Vepstas,
{\it On the Minkowski Measure}
arXiv preprint arXiv:0810.1265 (2008).


\end{thebibliography}
\end{document}